\documentclass[11pt]{article}
\usepackage{epsfig}
\usepackage{graphicx}
\usepackage[latin1]{inputenc}
\usepackage[T1]{fontenc}
\usepackage{mathptmx}
\usepackage{url}

\usepackage{enumitem}

\usepackage{amsmath}

\usepackage{pdfpages}
\newtheorem{Def}{Definition}
\newtheorem{Thm}[Def]{Theorem}
\newtheorem{Le}[Def]{Lemma}
\newtheorem{Ex}[Def]{Example}
\newtheorem{Rem}[Def]{Remark}
\newtheorem{Sit}[Def]{Situation}

\begin{document}

\title{On distinguishing special trees by their chromatic symmetric functions}
\author{Melanie Gerling}
\date{\today}

\maketitle

\begin{abstract}
In 1995, Stanley introduced the well-known chromatic symmetric function $X_{G}(x_{1},x_{2},\ldots)$ of a graph $G$. It is a sum of monomial symmetric functions such that for each vertex coloring of $G$ there is exactly one of these summands. The question, whether $X_{G}(x_1,x_{2},\ldots)$ distinguishes nonisomorphic trees with the same number of vertices, is still open in general. For special trees it has already been shown. In 2008, Martin, Morin and Wagner proved it for spiders and some caterpillars. We decompose a tree by separating the leafs and their neighbors and do the same to the remaining forest until there remains a forest with vertices of degree not greater than $1$. For nonisomorphic trees $G$ and $H$ with the same number of vertices and special properties concerning their number of leafs and leaf neighbors for each subgraph in their leaf decomposition we prove that the chromatic symmetric function distinguishes the graphs. Our idea is to find independent partitions of $G$ and $H$ with respectively a block of maximal cardinality for $G$ as well as for $H$. These cardinalities are different for graphs $G$ and $H$ with special properties of their leaf decompositions. Additionally, we give explicit formulas for the cardinality of such a maximal block of an independent partition of star connections and spiders.
\end{abstract}

\section{Introduction}

Let $M$ be a set and let $\mathcal{P}(M)=\left\{X | X\subseteq M\right\}$ denote the power set of $M$. For a natural number $k\geq1$ let $\mathcal{P}_{k}(M)=\left\{X\in\mathcal{P}(M) | \left|X\right|=k\right\}$. Because we will exclude certain types of graphs, a simplified definition of a graph suffices. We define a graph $G$ as a tuple $G=\left(V(G),E(G)\right)=\left(V,E\right)$ of two sets $V(G)=V$ and $E(G)=E$ with $E\subseteq\mathcal{P}_{2}(V)$. The members of $V$ are called \textbf{vertices} and the members of $E$ \textbf{edges} of $G$. This definition excludes the occurrence of \textbf{loops}, i. e. members of $E\cap \mathcal{P}_{1}(V)$, and \textbf{multiple edges}, which are irrelevant for our considerations. Two vertices $v_{1},v_{2}\in V$ are called \textbf{adjacent} to each other, if $\left\{v_{1},v_{2}\right\}\in E$. For a vertex $v\in V$ the number of all vertices that are adjacent to $v$ in $G$ is called the \textbf{degree} of $v$ in $G$ and denoted by \textbf{$deg_{G}(v)$}.\\
A graph $H=\left(V',E'\right)$ is called a \textbf{subgraph} of $G=\left(V,E\right)$, if $V'\subseteq V$ and $E'\subseteq E$. We also write $H\subseteq G$. A subgraph is an \textbf{induced subgraph} of $G$, if $E'=E\cap\mathcal{P}_{2}(V')$. For $X\subseteq V$ the graph $G[X]$ is the subgraph induced by the set $X$ and for $\left\{v_{1},\ldots,v_{r}\right\}\subseteq V$ we set $G-\left\{v_{1},\ldots,v_{r}\right\}=G[V\setminus\left\{v_{1},\ldots,v_{r}\right\}]$.\\
Let $G=(V,E)$ be a finite graph and $n=|V|$. A function $\kappa: V\rightarrow Y$ with a set of colors $Y=\left\{1,2,\ldots\right\}$ is called a \textbf{vertex coloring}. It is called a \textbf{proper vertex coloring}, if for all edges $\left\{v_{1},v_{2}\right\}$ the relation $\kappa(v_{1})\neq\kappa(v_{2})$ holds. The number of all proper vertex colorings for a graph $G$ in at most $y$ colors was introduced as a polynomial $P(G;y)$ by \cite{B} and was called the \textbf{chromatic polynomial} by \cite{BL}. Now, we turn to a similar function.

\begin{Def}(\cite{STAN}, p. 168)\\
Let $x_{1},x_{2},\ldots$ be (commuting) indeterminates and $V=\left\{v_{1},v_{2}\ldots,v_{n}\right\}$. Then the \textbf{symmetric function generalization of the chromatic polynomial} of $G$ is defined as
\begin{eqnarray*}
X_{G}(x_{1},x_{2},\ldots)=\sum_{\kappa}x_{\kappa(v_{1})}x_{\kappa(v_{2})}\ldots x_{\kappa(v_{n})}
\end{eqnarray*}
where $\kappa:V\rightarrow Y$ is a proper vertex coloring.
\label{SymStan}
\end{Def}

This function has the following connection with the chromatic polynomial.

\begin{Rem}(\cite{STAN}, p. 169)\\
For $r\in Y$ we have $X_{G}(1^{r})=P(G;r)$.
\end{Rem}

For the next Remark we need some definitions.

\begin{Def}
Let $G=\left(V,E\right)$ be a graph. We call a partition of $V$ a \textbf{partition of $G$}. \textbf{$\Pi(G)$} denotes the set of all partitions of $G$. Each member of a partition is called a \textbf{block}.
\end{Def}

We will only be interested in partitions which are defined as follows.

\begin{Def}
Non adjacent vertices and edges which don't have any vertex in common are called \textbf{independent of each other}. A set $X\subseteq V$ or $Y\subseteq E$ is said to be \textbf{independent} if its members are pairwise independent of each other.\\
An \textbf{independent partition} $\pi$ of $G=\left(V,E\right)$ denotes, as defined in \cite{DPT}, section 3, a partition of $G$, such that each block $A\in\pi$ is independent. We denote the set of all independent partitions by \textbf{$\Pi_{I}(G)$}.
\label{KnotPar.2}
\end{Def}

The next obersvation will be very useful to answer the question for special types of trees, whether two not isomorphic trees with the same number of vertices can be distinguished by their symmetric function generalization of the chromatic polynomial or not.

\begin{Rem}
We can identify a vertex coloring of a graph $G=(V,E)$ with an independent partition of $V$, such that the vertices of each block are colored in the same color and vertices from different blocks are colored in different colors. As Definition \ref{SymStan} shows, each summand $x_{\kappa(v_{1})}x_{\kappa(v_{2})}\ldots x_{\kappa(v_{n})}$ contains a product $x_{\kappa(v_{i})}^{k}$, if there are vertices $v_{1},\ldots,v_{k}\in V$ with $\kappa(v_{1})=\ldots=\kappa(v_{k})$, which means that they are in the same block. So, $k$ shows the cardinality of that block.
\label{IndParCol}
\end{Rem}

For special types of non-isomorphic trees with the same number of vertices it will be easy to show, that the cardinalities of their blocks with maximal cardinality are different from each other. So the chromatic symmetric functions of these trees differ of course, too.\\

We call a graph $P=\left(V,E\right)$ with vertex set $V=\left\{v_{0},\ldots,v_{k}\right\}$ and edge set $E=\left\{\left\{v_{0},v_{1}\right\},\left\{v_{1},v_{2}\right\},\ldots,\left\{v_{k-1},v_{k}\right\}\right\}$ with $|V|\geq1$ and pairwise distinct $v_{i}\textrm{ for }i=0,\ldots,k$ a \textbf{path} of length $k$ from $v_{0}$ to $v_{k}$.\\
A non-empty graph is called \textbf{connected}, if for all $u,w\in V$ there is a path from $u$ to $w$, i. e. a path $P$ with $u=v_{0}$ and $w=v_{k}$.\\
A graph $F$ without any cycles is called a \textbf{forest}, a connected forest $T$ is called a \textbf{tree} and a vertex $v\in V(F)$ of degree one is called a \textbf{leaf}.

\begin{Def}
Let $T=T_{1}=F_{1}$ be a tree, the vertices $v_{1},v_{2},\ldots,v_{b_{1}}\in V(T)$ its leaves and $w_{1},w_{2},\ldots w_{\eta_{1}}\in V(T)$ all vertices that are adjacent to a leaf. Of course, $b_{1}\geq\eta_{1}$ and $F_{1}-\left\{v_{1},\ldots,v_{b_{1}},w_{1},\ldots,w_{\eta_{1}}\right\}=F_{2}$ is a forest. For $F_{2}$ we do the same and continue in this manner. We get a subgraph relation $T_{1}=F_{1}\subseteq F_{2}\subseteq\ldots\subseteq F_{k}$, such that $F_{k}$ is a forest with vertex degrees $\leq1$. For each $i\in\left\{1,\ldots,k-1\right\}$ there are leaves $v_{i_{1}},v_{i_{2}},\ldots,v_{i_{b_{i}-1}}\in V(F_{i})$ and all possible vertices $w_{i_{1}},w_{i_{2}},\ldots,w_{i_{\eta_{i}-1}}\in V(F_{i})$, that are adjacent to a leaf. If in the forest $F_{k}$ all vertices have degree $0$, we have $b_{1}\geq \eta_{1}\geq b_{2}\geq \eta_{2}\geq\ldots\geq b_{k}\geq \eta_{k}$ with $\eta_{k}=0$. But if there are $\alpha$ vertices $v_{kb_{k}1},\ldots,v_{kb_{k}\alpha}$ of degree $1$ with $1\leq\alpha\leq|V(F_{k})|$, we set $b_{k}=|V(F_{k})|-\frac{\alpha}{2}$ and $\eta_{k}=\frac{\alpha}{2}$. Note, that $\alpha$ is always even. Then $b_{1}\geq\eta_{1}\geq b_{2}\geq\eta_{2}\geq\ldots \geq b_{k}\geq\eta_{k}$ holds. For a tree $T$ we call such a relation a \textbf{leaf decomposition} and denote it by $(T,k)$.
\label{def.comp}
\end{Def}

\begin{Rem}
For every tree $T$ there is exactly one $k\in\mathbf{N}$, $k\geq0$, such that $(T,k)$ is a leaf decomposition of $T$ as in Definition \ref{def.comp}.\\
Each vertice $v\in V(T)$ is of the form $v_{ij_{1}}$ or $w_{ij_{2}}$ with $i\in\left\{1,\ldots,k\right\}$, $j_{1}\in\left\{1,\ldots,b_{i}\right\}$ and $j_{2}\in\left\{1,\ldots,\eta_{i}\right\}$. So, $v$ is either a leaf or a neighbor of a leaf in the associated graph $F_{i}$.
\end{Rem}

Note, that two non-isomorphic trees are not determined by their tree decomposition as the following example shows.

\begin{Ex}
In the following picture there are two non-isomorphic trees $T_{1}$ and $T_{2}$ with the same leaf decomposition. For $T_{1}$ we define the leaf decomposition $b_{11}\geq \eta_{11}\geq b_{12}$ and for $T_{2}$ the leaf decomposition $b_{21}\geq \eta_{21}\geq b_{22}$. We have $b_{11}=b_{21}=6$, $\eta_{11}=\eta_{21}=3$ and $b_{12}=b_{22}=1$.
\begin{figure}[ht]
\hspace*{3cm}\includegraphics[width=10cm]{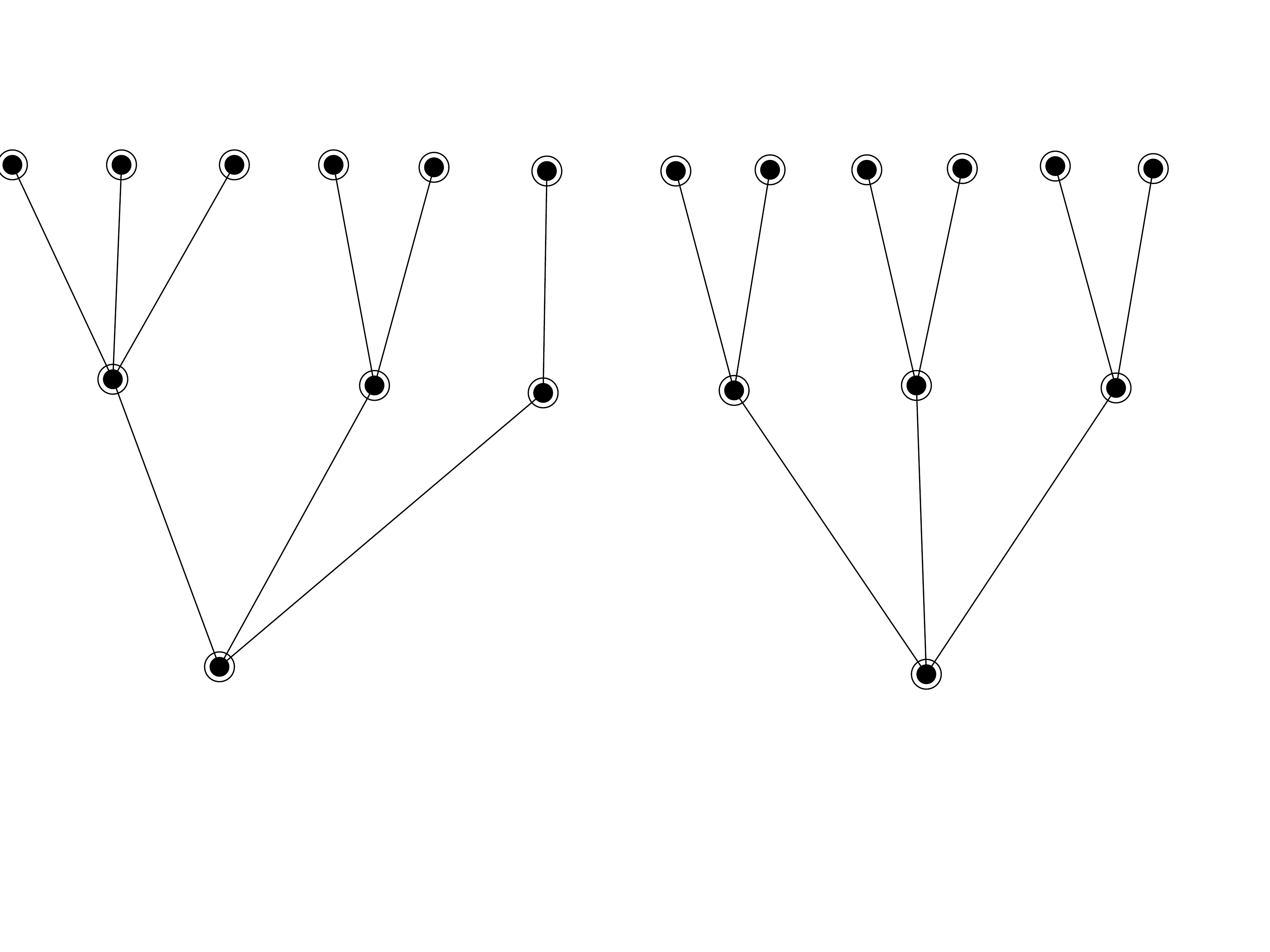}
\caption{$T_{1}$ and $T_{2}$}
\end{figure}
\label{compex1}
\end{Ex}

\section{The number of leaves and further properties}

\begin{Thm}
Let $n\geq3$. For a tree $T$ with $|V(T)|=n$ there is a maximum $|B|$ with $B\in\sigma$, $\sigma\in\Pi_{I}(T)$, such that for all leaves $v\in V(T)$ we have $v\in B$.
\label{MaxBlock1}
\end{Thm}

\textbf{Proof}:\\
For all leaves $v\in V(T)$ there is a vertex $w$ with $N_{T}(v)=\left\{w\right\}$, i. e. of course $|N_{T}(v)|=1$. From all such pairs $\left\{v,w\right\}$ we may choose at most one vertex for our block $B$ because of its independence and because of the maximality of $|B|$ we have to chose at least one. We choose always the leaf $v$ because the number of neighbors of leafs is equal to or smaller then the number of leafs. So we have $v\in B$ for all leaves $v\in V(T)$.\\
q.e.d.\\
\\
The number of leaves can be an important hint to decide, whether two trees can be distinguished by their symmetric function generalization of the chromatic polynomial as the next theorem shows.

\begin{Thm}
Let $T_{1}$ and $T_{2}$ be trees with $|V(T_{1})|=|V(T_{2})|=n\geq4$. The number of leaves of $T_{i}$ for $i\in\left\{1,2\right\}$ is denoted by $b_{i}$ and $b_{1}>b_{2}$. Let $\eta_{i}$ be the number of all vertices $w$ with $v\in N_{T_{i}}(w)$ for a leaf $v\in V(T_{i})$ and let $\rho_{i}=n-\eta_{i}-b_{i}$. Let $V(\rho_{i})$ the set of vertices which contribute to $\rho_{i}$ and let $T[V(\rho_{i})]$ be a path. Furthermore, we assume one of the following cases:
\begin{enumerate}
\item
$\rho_{1}=\rho_{2}$,
\item
$\rho_{1}>\rho_{2}$,
\item
$\rho_{1}<\rho_{2}$ and $b_{1}-b_{2}>\left\lceil \frac{\rho_{2}-\rho_{1}}{2}\right\rceil$,
\item
$\rho_{1}<\rho_{2}$ and $b_{1}-b_{2}>\left\lceil\frac{\rho_{1}-\rho_{2}}{k}\right\rceil$ for all $k\in\mathrm{N}$ with $k\neq2$ and\\
$\frac{k}{k-2}\leq\rho_{2}-\rho_{1}<k(b_{1}-b_{2})$.
\end{enumerate}
Then $X_{T_{1}}(x_{1}, x_{2},\ldots)\neq X_{T_{2}}(x_{1}, x_{2},\ldots)$.
\label{DisLeaves1}
\end{Thm}

\textbf{Proof}:\\
Because of Theorem \ref{MaxBlock1} we compare independent partitions of $T_{1}$ and $T_{2}$, such that there is respectively one block $B_{i}$ with $|B_{i}|$ maximum and $v_{i}\in B_{i}$ for all leaves $v_{i}\in V(T_{i})$. We set 
\begin{eqnarray*}
M_{i}&=&max\left\{|B_{i}| | B_{i}\in\sigma_{i} \wedge \sigma{i}\in\Pi_{I}(T_{i})\right\}\;.
\end{eqnarray*}
For each case we will show $M_{1}>M_{2}$. Note, that we can choose from the vertices $\rho_{i}$ exactly $\lceil\frac{\rho_{i}}{2}\rceil$ further vertices for the block $B_{i}$ because we did not choose any neighbor of a leaf and because of the path structure of $T[V(\rho_{i})]$.
\begin{enumerate}
\item
\begin{enumerate}
\item
Case $\rho_{1}$ is even. Because of $b_{1}>b_{2}$ and $\rho_{1}=\rho_{2}$ we have
\begin{eqnarray*}
M_{1}=b_{1}+\frac{\rho_{1}}{2}&>&b_{2}+\frac{\rho_{2}}{2}=M_{2}\;.
\end{eqnarray*}
\item
Case $\rho_{1}$ is odd. Like in the first case we get
\begin{eqnarray*}
M_{1}=b_{1}+\frac{\rho_{1}+1}{2}&>&b_{2}+\frac{\rho_{2}+1}{2}=M_{2}\:.
\end{eqnarray*}
\end{enumerate}
\item
\begin{enumerate}
\item
Case $\rho_{1}$ and $\rho_{2}$ are even. Because of $b_{1}>b_{2}$ and $\rho_{1}>\rho_{2}$ it leads to
\begin{eqnarray*}
M_{1}=b_{1}+\frac{\rho_{1}}{2}&>&b_{2}+\frac{\rho_{2}}{2}=M_{2}\;.
\end{eqnarray*}
\item
Case $\rho_{1}$ and $\rho_{2}$ are odd. Like above the following holds:
\begin{eqnarray*}
M_{1}=b_{1}+\frac{\rho_{1}+1}{2}&>&b_{2}+\frac{\rho_{2}+1}{2}=M_{2}\;.
\end{eqnarray*}
\item
Case $\rho_{1}$ is odd and $\rho_{2}$ is even. Similar as above we get
\begin{eqnarray*}
M_{1}=b_{1}+\frac{\rho_{1}+1}{2}&>&b_{2}+\frac{\rho_{2}}{2}=M_{2}\;.
\end{eqnarray*}
\item
Case $\rho_{1}$ is even and $\rho_{2}$ is odd. Like above and because of $\rho_{1}\geq\rho_{2}+1$ it leads to
\begin{eqnarray*}
M_{1}=b_{1}+\frac{\rho_{1}}{2}&>&b_{2}+\frac{\rho_{2}+1}{2}=M_{2}\;.
\end{eqnarray*}
\end{enumerate}
\item
\begin{enumerate}
\item
Case $\rho_{1}$ and $\rho_{2}$ are even. Hence,
\begin{eqnarray*}
b_{1}-b_{2}&>&\left\lceil\frac{\rho_{2}-\rho_{1}}{2}\right\rceil\\
M_{1}=b_{1}+\frac{\rho_{1}}{2}&>&b_{2}+\frac{\rho_{2}}{2}=M_{2}\;.
\end{eqnarray*}
\item
Case $\rho_{1}$ and $\rho_{2}$ are odd. Similary we get
\begin{eqnarray*}
b_{1}-b_{2}&>&\left\lceil\frac{\rho_{2}-\rho_{1}}{2}\right\rceil\\
M_{1}=b_{1}+\frac{\rho_{1}+1}{2}&>&b_{2}+\frac{\rho_{2}+1}{2}=M_{2}\;.
\end{eqnarray*}
\item
Case $\rho_{1}$ is odd and $\rho_{2}$ is even. Like above we have
\begin{eqnarray*}
b_{1}-b_{2}&>&\frac{\rho_{2}-\rho_{1}}{2}-\frac{1}{2}\\
M_{1}=b_{1}+\frac{\rho_{1}+1}{2}&>&b_{2}+\frac{\rho_{2}}{2}=M_{2}\;.
\end{eqnarray*}
\item
Case $\rho_{1}$ is even and $\rho_{2}$ is odd. Thus,
\begin{eqnarray*}
b_{1}-b_{2}&&>\left\lceil\frac{\rho_{2}-\rho_{1}}{2}\right\rceil\\
b_{1}-b_{2}&>&\frac{\rho_{2}-\rho_{1}}{2}+\frac{1}{2}\\
M_{1}=b_{1}+\frac{\rho_{1}}{2}&>&b_{2}+\frac{\rho_{2}+1}{2}=M_{2}\;.
\end{eqnarray*}
\end{enumerate}
\item
First we show that $\frac{\rho_{2}-\rho_{1}}{k}\geq\frac{\rho_{2}-\rho_{1}}{2}+\frac{1}{2}$ holds. As assumed, we have
\begin{eqnarray*}
\rho_{2}-\rho_{1}&\geq&\frac{k}{k-2}\\
\Leftrightarrow(k-2)(\rho_{2}-\rho_{1})-k&\geq&0\\
\Leftrightarrow(2-k)\rho_{2}+(k-2)\rho_{1}-k&\geq&0\\
\Leftrightarrow2\rho_{2}-2\rho_{1}&\geq&k\rho_{2}-k\rho_{1}+k\\
\Leftrightarrow\frac{\rho_{2}-\rho_{1}}{k}&\geq&\frac{\rho_{2}-\rho_{1}}{2}+\frac{1}{2}\;.
\end{eqnarray*}
\begin{enumerate}
\item
Case $\rho_{1}$ and $\rho_{2}$ are even. Then we get
\begin{eqnarray*}
b_{1}-b_{2}&\geq&\left\lceil\frac{\rho_{2}-\rho_{1}}{k}\right\rceil\geq\left\lceil\frac{\rho_{2}-\rho_{1}}{2}\right\rceil\\
M_{1}=b_{1}+\frac{\rho_{1}}{2}&\geq&b_{2}+\frac{\rho_{2}}{2}=M_{2}\;.
\end{eqnarray*}
\item
Case $\rho_{1}$ and $\rho_{2}$ are odd. It leads to
\begin{eqnarray*}
b_{1}-b_{2}&\geq&\left\lceil\frac{\rho_{2}-\rho_{1}}{k}\right\rceil\geq\left\lceil\frac{\rho_{2}-\rho_{1}}{2}\right\rceil\\
b_{1}-b_{2}&\geq&\frac{\rho_{2}+1}{2}-\frac{\rho_{1}+1}{2}\\
M_{1}=b_{1}+\frac{\rho_{1}+1}{2}&\geq&b_{2}+\frac{\rho_{2}+1}{2}=M_{2}\;.
\end{eqnarray*}
\item
Case $\rho_{1}$ is odd and $\rho_{2}$ is even. Similary as above we have
\begin{eqnarray*}
b_{1}-b_{2}&\geq&\left\lceil\frac{\rho_{2}-\rho_{1}}{k}\right\rceil\geq\left\lceil\frac{\rho_{2}-\rho_{1}}{2}\right\rceil\\
b_{1}-b_{2}&\geq&\frac{\rho_{2}-\rho_{1}}{2}-\frac{1}{2}\\
M_{1}=b_{1}+\frac{\rho_{1}+1}{2}&\geq&b_{2}+\frac{\rho_{2}}{2}=M_{2}\;.
\end{eqnarray*}
\item
Case $\rho_{1}$ is even and $\rho_{2}$ is odd. Hence,
\begin{eqnarray*}
b_{1}-b_{2}&\geq&\left\lceil\frac{\rho_{2}-\rho_{1}}{k}\right\rceil\geq\left\lceil\frac{\rho_{2}-\rho_{1}}{2}\right\rceil\\
b_{1}-b_{2}&\geq&\frac{\rho_{2}-\rho_{1}}{2}+\frac{1}{2}\\
M_{1}=b_{1}+\frac{\rho_{1}}{2}&\geq&b_{2}+\frac{\rho_{2}+1}{2}=M_{2}\;.
\end{eqnarray*}
\end{enumerate}
\end{enumerate}
q.e.d.\\
\\
We show an example for Theorem \ref{DisLeaves1}. It contains two graphs which are special cases of the graphs claimed in Theorem \ref{DisLeaves1}. This special graph type was already introduced by \cite{MMW}. They already proved that so-called caterpillars can be distinguished by their chromatic symmetric functions.

\begin{Ex}
Let $n=4$. We consider the following two graphs. On the left side there is the graph $S_{4}$ and on the right side we see the graph $P_{4}$. We have $|V(S_{4})|=|V(P_{4})|=4$.
\begin{figure}[ht]
\hspace*{3cm}\includegraphics[width=10cm]{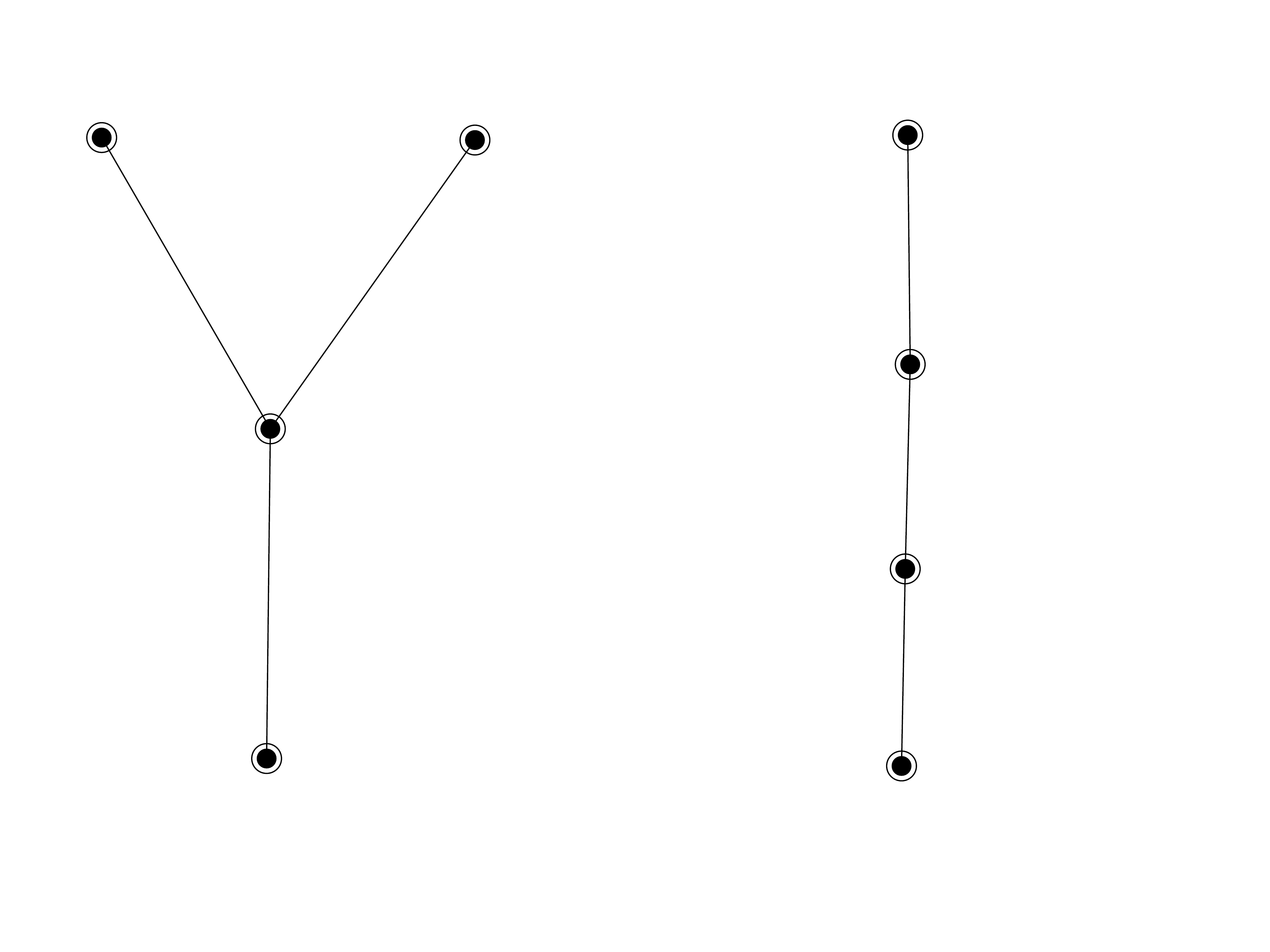}
\caption{$S_{4}$ and $P_{4}$}
\end{figure}
For the graph $S_{4}$ we have $b_{1}=3$, $\eta_{1}=1$ and $\rho_{1}=0$, while we have $b_{2}=2$, $\eta_{2}=2$ and $\rho_{2}=0$ for the graph $P_{4}$. Because of $b_{1}>b_{2}$ and $\rho_{1}=\rho_{2}$ Theorem \ref{DisLeaves1} leads to $3=M_{1}>M_{2}=2$ as we can see easily.
\label{Bild.Beispiel.1}
\end{Ex}

Now, we will pay attention to a special situation which is interesting for our further considerations.

\begin{Sit}
Let $T_{11}$ and $T_{21}$ be trees with $|V(T_{11})|=|V(T_{21})|=n$ and the leaf decompositions $(T_{11},k_{1})$ and $(T_{21},k_{2})$. We set $r=\textrm{ max }\left\{k_{1},k_{2}\right\}$. Let $T_{11}\subseteq F_{12}\subseteq F_{13}\subseteq\ldots\subseteq F_{1k_{1}}$, and for $i\in\left\{1,\ldots,k_{1}\right\}$ there exist respectively $b_{1i}$ leafs and $\eta_{1i}$ neighbor vertices like in Definition \ref{def.comp}. Analogously, let $T_{21}\subseteq F_{22}\subseteq F_{23}\subseteq\ldots\subseteq F_{2k_{2}}$ with $b_{2j}$ and $\eta_{2j}$ for $j\in\left\{1,\ldots,k_{2}\right\}$ like in Definition \ref{def.comp}.\\
We set $(T_{11},k_{1})=(T_{11},r)$ and $(T_{21},k_{2})=(T_{21},r)$, where $b_{1i}$ and $\eta_{1i}$ or $b_{2j}$ and $\eta_{2j}$ become zero for $k_{1}\neq k_{2}$ from a certain $i$ or respectively $j$ on.
\label{sit1}
\end{Sit}

We turn to some special cases. The next theorem shows the first one.

\begin{Thm}
We regard Situation \ref{sit1} for trees $T_{1}$ and $T_{2}$ and assume $\eta_{1i}\leq\eta_{2i}$ as well as $b_{1i}\geq b_{2i}$ for all $i\in\left\{1,\ldots,r\right\}$. We except the case that there is equality everywhere. Then $X_{T_{1}}(x_{1},\ldots)\neq X_{T_{2}}(x_{1},\ldots)$.
\label{satz1}
\end{Thm}

\textbf{Proof}:\\
By Theorem \ref{MaxBlock1}, for $k=1,2$ there is a partition $\sigma_{k}\in\Pi_{I}(T_{k})$ with $B_{k}\in\sigma_{k}$ and $|B_{k}|$ maximal, such that for all leaves $v\in T_{k}$ we have $v\in B_{k}$. Let $M_{k}$ as in the proof of Theorem \ref{DisLeaves1}. Then we have $b_{ki}\in B_{k}$, $\eta_{ki}\notin B_{k}$ for all $i\in\left\{1,\ldots,r\right\}$. Because we have at least one strict inequality the relation $M_{1}>M_{2}$ holds because of $|V(T_{1})|=|V(T_{2})|=n$ and so $X_{T_{1}}(x_{1},\ldots)\neq X_{T_{2}}(x_{1},\ldots)$.\\
q.e.d.\\

\begin{Ex}
Let $n=8$. We consider the following two graphs. On the left side there is a tree $T_{1}$ and on the right side we see another tree $T_{2}$. We have $|V(T_{1})|=|V(T_{2})|=8$.
\begin{figure}[ht]
\hspace*{3cm}\includegraphics[width=10cm]{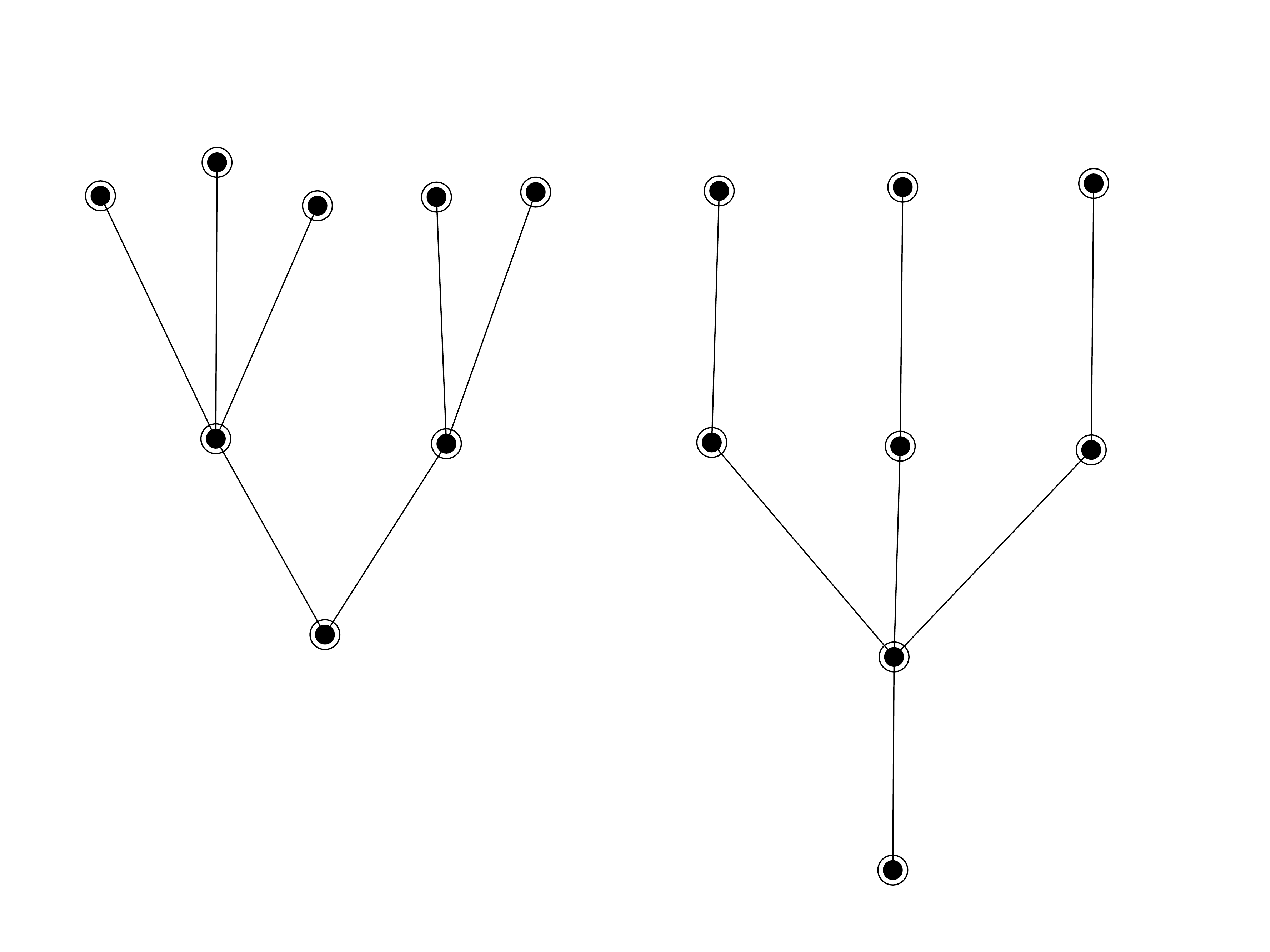}
\caption{$T_{1}$ and $T_{2}$}
\end{figure}
We get the following relations:
\begin{eqnarray*}
b_{11}=5&>&b_{21}=4\\
\eta_{11}=2&<&\eta_{21}=4\\
b_{12}=1&>&b_{22}=0\\
\eta_{12}=0&=&\eta_{22}=0\;,
\end{eqnarray*}
such that Theorem \ref{satz1} holds. Indeed,
\begin{eqnarray*}
|M_{1}|=5+1=6&>&|M_{2}|=4+0=4.
\end{eqnarray*}
\label{Bild.Beispiel.Satz.2}
\end{Ex}

In the following Theorem we introduce another case.

\begin{Thm}
We consider again Situation \ref{sit1}. Let $i_{1},i_{2}\ldots,i_{s}\in\left\{1,2,\ldots,r\right\}$, $1\leq s\leq r-1$ and $b_{1j}>b_{2j}$ for all $j\in\left\{1,2,\ldots,r\right\}\setminus\left\{i_{1},i_{2},\ldots,i_{s}\right\}$ and $b_{1i_{t}}\leq b_{2i_{t}}$ for all $t\in\left\{1,2,\ldots,s\right\}$. Furthermore,
\begin{eqnarray}
\sum_{\alpha=i_{1}}^{i_{s}}(b_{2\alpha}-b_{1\alpha})&<&\sum_{j\in\left\{1,\ldots,r\right\}\setminus\left\{i_{1},\ldots,i_{s}\right\}}(b_{1j}-b_{2j})
\label{eqn:ineq}
\end{eqnarray}
holds. Then we have $X_{T_{1}}(x_{1},\ldots)\neq X_{T_{2}}(x_{1},\ldots)$.
\label{satz2}
\end{Thm}

\textbf{Proof}:\\
As in the proof of Theorem \ref{satz1}, we look at $\sigma_{k}\in\Pi_{I}(T_{k})$ for $k=1,2$ with $B_{k}\in\sigma_{k}$ and $|B_{k}|$ maximal, such that for all leaves $v\in T_{k}$ we have $v\in B_{k}$. Again, $b_{ki}\in B_{k}$, $\eta_{ki}\notin B_{k}$ for all $i\in \left\{1,\ldots,r\right\}$. Because of inequality \ref{eqn:ineq} we have $M_{1}>M_{2}$ and hence $X_{T_{1}}(x_{1},\ldots)\neq X_{T_{2}}(x_{1},\ldots)$.
\\
q.e.d.\\

\section{Maximal blocks of independent partitions of star connections}

In this section we will be concerned with maximal blocks of independent partitions of a special type of graph. It is constructed from stars, which are renownedly defined as follows.

\begin{Def}
A \textbf{star $S_{n}$} for a natural number $n\geq 2$ is a graph with vertex set $V=\left\{v_{1},\ldots,v_{n}\right\}$ and edge set $E=\left\{\left\{v_{1},v_{2}\right\},\left\{v_{1},v_{3}\right\},\ldots,\left\{v_{1},v_{n}\right\}\right\}$.\\
The vertex $v_{1}$ is called the \textbf{center} of the star.
\label{Star}
\end{Def}

Now, we introduce two operations. The first one is a graph operation and the second one a set operation.\\
Let $G_{1}=(V_{1},E_{1})$ and $G_{2}=(V_{2},E_{2})$ be two graphs. Then we can create a new graph $G_{1}\cup G_{2}=\left(V_{1}\cup V_{2}, E_{1}\cup E_{2}\right)$. Note, that the sets $V_{1}$ and $V_{2}$ (and thus of course possibly $E_{1}$ and $E_{2}$) are not necessarily disjoint.\\
An intersection $\bigcap_{i\in\left\{1,\ldots,s\right\}}A_{i}$ of sets, such that the $i\in\left\{1,\ldots,s\right\}$ are pairwise distinct from each other, will be denoted by $\bigcap_{i\in\left\{1,\ldots,s\right\}}^{\bullet}A_{i}$. For a special type of graph we give the following Definition.

\begin{Def}
Let $S_{n_{1}},\ldots,S_{n_{r}}$ be stars with $|V(S_{n_{i}})|=n_{i}$ vertices for $i=1,\ldots,r$ with $r\geq2$ and $n_{i}\geq3$. We define a tree $T=(V,E)=\bigcup_{i=1}^{r}S_{n_{i}}$ with the following properties: Let $C$ be the set of all centers of stars in $T$ and 
\begin{eqnarray*}
&&\left\{v_{1},\ldots,v_{t}\right\}=\\
&&\left\{v\in V\setminus C\textrm{  } | \textrm{  } \exists \textrm{  } S_{n_{i_{1}}},\ldots,S_{n_{i_{s}}} \textrm{ with } i_{1},\ldots,i_{s}\in\left\{1,\ldots,r\right\}, \textrm{  } \left\{v\right\}=\bigcap_{l\in\left\{i_{1},\ldots,i_{s}\right\}}^{\bullet}S_{n_{l}}\right\}\;.
\end{eqnarray*}
Then we assume $|S_{n_{i}}|\cap|S_{n_{j}}|\leq1$ for all pairwise distinct $i,j\in\left\{1,\ldots,r\right\}$ and that no center belongs to more than one star as a center or another vertex. Hence, for each $v\in V\setminus (C\cup\left\{v_{1},\ldots,v_{t}\right\})$ we have $deg_{G}(v)=1$ and for $v\in\left\{v_{1},\ldots,v_{t}\right\}$ with $\left\{v\right\}=\bigcap_{l\in\left\{i_{1},\ldots,i_{s}\right\}}^{\bullet}S_{n_{l}}$ we have $deg_{G}(v)=s$. Then we call the graph $T$ a \textbf{star connection} and denote it by $(S_{n_{1}},\ldots,S_{n_{r}},v_{1},\ldots,v_{t})$.
\label{StarUnion1}
\end{Def}

Before we continue with the symmetric function generalization of the chromatic polynomial, we observe the following property of a star connection.

\begin{Thm}
Like in the proof of Theorem \ref{DisLeaves1}, we define 
\begin{eqnarray*}
M&=&max\left\{|B_{j}| | B_{j}\in\sigma_{j}\wedge\sigma_{j}\in\Pi_{I}(T)\right\}\;.
\end{eqnarray*}
In the situation of Definition \ref{StarUnion1} we have
\begin{eqnarray*}
M&=&\sum_{k=1}^{r}(n_{k}-1)-\sum_{i=1}^{t}(deg_{T}(v_{i})-1)\;.
\end{eqnarray*}
\label{StarUnion2}
\end{Thm}

\textbf{Proof}:\\
We will prove the Theorem by complete induction.\\
For $r=2$ the graph $T$ is a star connection $S_{n_{1}}\cup S_{n_{2}}$ with exactly one common vertex $v_{1}$, that is not a center of a star. Thus, $deg_{G}(v_{1})=2$. The greatest set of pairwise independent vertices is of course the set $A$ of all vertices without the two centers. Additionally, we may count $v_{1}$ only one time. Hence, for $|A|=M$ we have
\begin{eqnarray*}
M&=&n_{1}+n_{2}-3\\
&=&n_{1}+n_{2}-2-(2-1)\\
&=&\sum_{k=1}^{2}(n_{k}-1)-\sum_{i=1}^{1}(deg_{T}(v_{i})-1)\;.
\end{eqnarray*}
For the induction step we assume that the claim is true for an $r\geq2$. Let $T=\bigcup_{i=1}^{r}S_{n_{i}}$. We want to connect $T$ with a star $S_{n_{r+1}}$ to a new star connection $H$. Then there is a leaf $w_{1}\in V(S_{n_{r+1}})$ which is also a leaf of stars $S_{n_{i_{1}}},\ldots,S_{n_{i_{s}}}$ in $T$ with $i_{1},\ldots,i_{s}\in\left\{1,\ldots,r\right\}$ and $1\leq s\leq r$. We consider the maximal set $A_{T}$ of pairwise independent vertices in the graph $T$. Because $w_{1}$ is not a center we have $w_{1}\in A_{T}$. By induction, for $|A_{T}|=M_{T}$ we get
\begin{eqnarray*}
M_{T}&=&\sum_{k=1}^{r}(n_{k}-1)-\sum_{i=1}^{t}(deg_{T}(v_{i})-1)\;.
\end{eqnarray*}
All the other leaves $w_{2},\ldots,w_{n_{r+1}}$ of $S_{n_{r+1}}$ are not leaves of other stars. Hence, there are $n_{r+1}-2$ new vertices in $H$ that are contained in the maximal set $A_{H}$ of pairwise independent vertices in $H$. We have $deg_{H}(w_{1})=deg_{T}(w_{1})+1$.\\
Case 1: $w_{1}$ is identified in $H$ with a leaf of $T$ that is not a connection vertex in $T$. We set $w_{1}=v_{t+1}$. Then $deg_{T}(v_{i})=deg_{H}(v_{i})$ for all $i\in\left\{1,\ldots,t\right\}$ while $deg_{T}(v_{t+1})=1$ and $deg_{H}(v_{t+1})=2$, such that
\begin{eqnarray*}
M_{H}&=&M_{T}+(n_{r+1}-2)\\
&=&\sum_{k=1}^{r}(n_{k}-1)-\sum_{i=1}^{t}(deg_{T}(v_{i})-1)+(n_{r+1}-2)\\
&=&\sum_{k=1}^{r+1}(n_{k}-1)-\sum_{i=1}^{t}(deg_{T}(v_{i})-1)-1\\
&=&\sum_{k=1}^{r+1}(n_{k}-1)-\sum_{i=1}^{t}(deg_{H}(v_{i})-1)-(deg_{H}(v_{t+1})-1)\\
&=&\sum_{k=1}^{r+1}(n_{k}-1)-\sum_{i=1}^{t+1}(deg_{H}(v_{i})-1)\;.
\end{eqnarray*}
Case 2: $w_{1}$ is identified in $H$ with a vertex that is not a center and not a leaf of $T$. Then that vertex is already a connection vertex in $T$, say $v_{t}$. As in case 1, $deg_{T}(v_{i})=deg_{H}(v_{i})$ for all $i\in\left\{1,\ldots,t-1\right\}$ and $deg_{T}(v_{t})+1=deg_{H}(v_{t})$, such that
\begin{eqnarray*}
M_{H}&=&M_{T}+(n_{r+1}-2)\\
&=&\sum_{k=1}^{r}(n_{k}-1)-\sum_{i=1}^{t}(deg_{T}(v_{i})-1)+(n_{r+1}-2)\\
&=&\sum_{k=1}^{r+1}(n_{k}-1)-\sum_{i=1}^{t}(deg_{T}(v_{i})-1)-1\\
&=&\sum_{k=1}^{r+1}(n_{k}-1)-\sum_{i=1}^{t}(deg_{H}(v_{i})-1)\;.
\end{eqnarray*}
\\
q.e.d.\\
\\
Now we will give an example for Theorem \ref{StarUnion2}.

\begin{Ex}
We look at the graph $T=S_{4}\cup S_{5}\cup S_{3}\cup S_{4}$ with $|V(T)|=13$ vertices, the number of connection vertices $t=3$ and the number of centers $C=4$ which is depicted below.
\begin{figure}[ht]
\hspace*{3cm}\includegraphics[width=10cm]{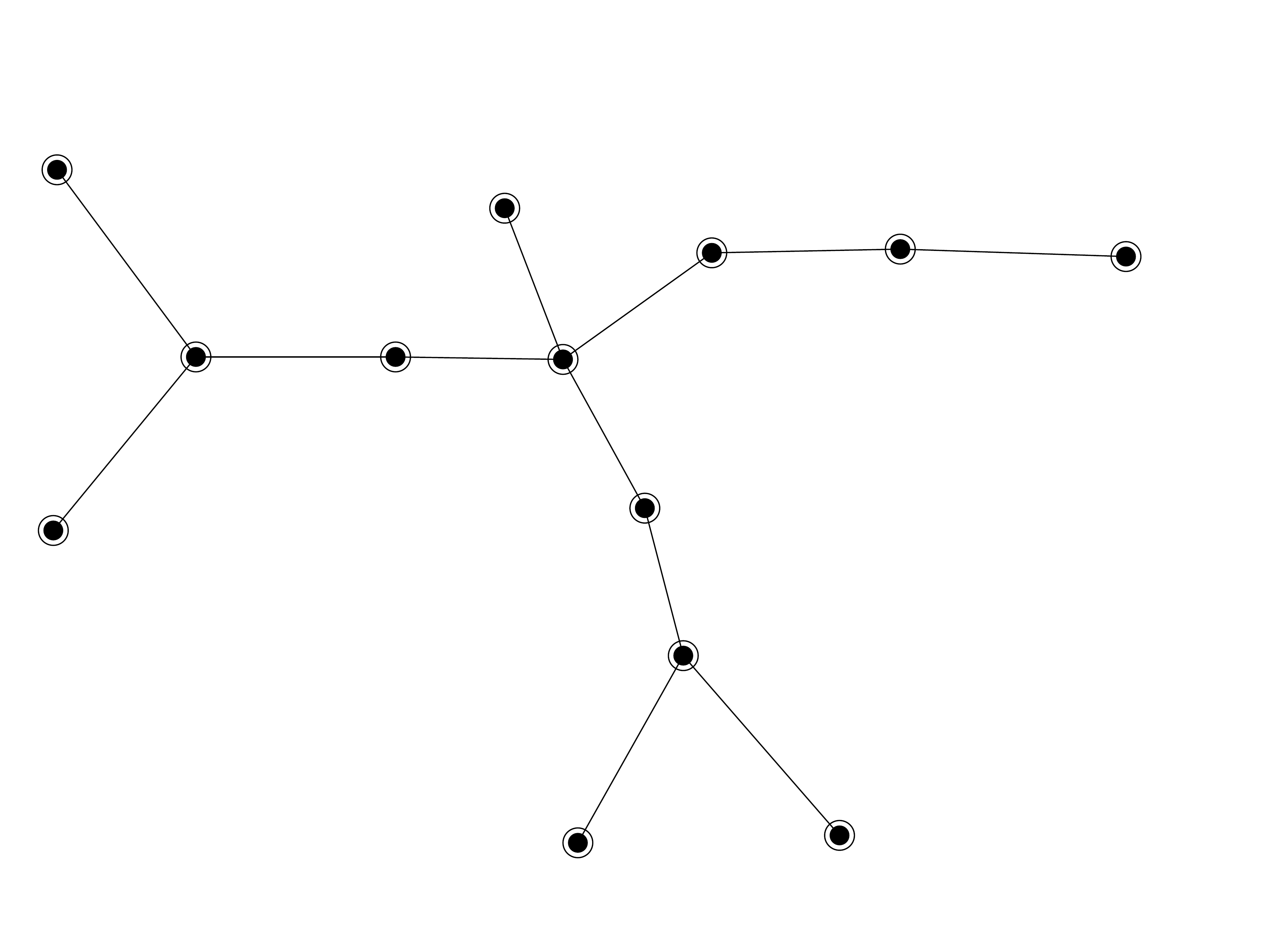}
\caption{$T$}
\label{Bild.Beispiel.StarUn3}
\end{figure}
In this graph we have $r=4$, and for the connection vertices $v_{1},\ldots,v_{3}$ we can see that $deg_{T}(v_{1})=deg_{T}(v_{2})=deg_{T}(v_{3})=2$. This leads to
\begin{eqnarray*}
M&=&\sum_{k=1}^{r}(n_{k}-1)-\sum_{i=1}^{t}(deg_{T}(v_{i})-1)\\
&=&(3+4+2+3)-(1+1+1)\\
&=&9\;.
\end{eqnarray*}
\label{StarUn3}
\end{Ex}

Now we can easily compare the symmetric function generalizations of the chromatic polynomial of two non-isomorphic star connections with the same number of vertices and get the following result. 

\begin{Thm}
Let $T_{1}=(S_{n_{1}},\ldots,S_{n_{r}},v_{1},\ldots,v_{t})$ and $T_{2}=(S_{m_{1}},\ldots,S_{m_{s}},w_{1},\ldots,w_{u})$ two star connections with $|V(T_{1})|=|V(T_{2})|$ and $r<s$. Then $X_{T_{1}}(x_{1},\ldots)\neq X_{T_{2}}(x_{1},\ldots)$.
\label{StarUn1}
\end{Thm}

To prove Theorem \ref{StarUn1} we need a Lemma.

\begin{Le}
In the situation of Definition \ref{StarUnion1} the two following relations hold.\\
\begin{itemize}
\item[a)] $|V(T)|=(\sum_{k=1}^{r}n_{k})-(r-1)$.
\item[b)]  $\sum_{i=1}^{t}(deg_{T}(v_{i})-1)=r-1$.
\end{itemize}
\label{StarUn2}
\end{Le}

\textbf{Proof}:\\
\begin{itemize}
\item[a)] Because of $V(T)=V(\bigcup_{i=1}^{r}S_{n_{i}})$ we have just to sum up all $n_{i}$ for all $i\in\left\{1,\ldots,r\right\}$ and to substract all $r-1$ multiple countings of vertices that connect different stars.
\item[b)] We show $\sum_{i=1}^{t}deg_{T}(v_{i})=t+r-1$.\\
For each vertex $v_{i}$  with $i\in\left\{1\ldots,t\right\}$ we have already the contribution $1=deg_{S_{n_{i}}}(v_{i})$, altogether $t$. If we construct $T$ step by step by the series $S_{n_{1}}, S_{n_{1}}\cup S_{n_{2}}, \ldots, S_{n_{1}}\cup\ldots\cup S_{n_{r}}=T$, we have $r-1$ steps. Each of them increases the degree of one of the vertices $v_{1},\ldots,v_{t}$ by $1$. Thus, $\sum_{i=1}^{t}deg_{T}(v_{i})=t+r-1$.
\end{itemize}
q.e.d.\\
\\
Now, we are able to prove Theorem \ref{StarUn1}.\\
\\
\textbf{Proof of Theorem \ref{StarUn1}}:\\
We show $M_{1}>M_{2}$.\\
By Theorem \ref{StarUnion2} and Lemma \ref{StarUn2} (b) we have
\begin{eqnarray*}
M_{1}&=&\sum_{k=1}^{r}(n_{k}-1)-\sum_{i=1}^{t}(deg_{T_{1}}(v_{i})-1)\\
&=&\sum_{k=1}^{r}(n_{k}-1)-(r-1)\\
&=&(\sum_{k=1}^{r}n_{k})-2r+1
\end{eqnarray*}
and
\begin{eqnarray*}
M_{2}&=&\sum_{k=1}^{s}(m_{k}-1)-\sum_{i=1}^{u}(deg_{T_{2}}(w_{i})-1)\\
&=&\sum_{k=1}^{s}(m_{k}-1)-(s-1)\\
&=&(\sum_{k=1}^{s}m_{k})-2s+1\;.
\end{eqnarray*}
By assumption and Lemma \ref{StarUn2} (a) we have furthermore 
\begin{eqnarray*}
|V(T_{1})|&=&(\sum_{k=1}^{r}n_{k})-(r-1)=|V(T_{2})|=(\sum_{k=1}^{s}m_{k})-(s-1)\;.
\end{eqnarray*}
Because of $r<s$
\begin{eqnarray*}
M_{1}=(\sum_{k=1}^{r}n_{k})-2r+1&>&(\sum_{k=1}^{s}m_{k})-2s+1=M_{2}
\end{eqnarray*}
follows.\\
q.e.d.

\section{Maximal blocks of independent partitions of spiders}

In this section we will turn our attention to so-called spiders which were introduced by \cite{MMW}. A \textbf{spider} is a tree which has exactly one vertex with degree greater than $2$.

\begin{Sit}
Let $S=(V,E)$ be a spider and $v$ the vertex with $r=deg_{S}(v)>2$. We denote the set of leafs of $S$ by $\left\{l_{1},\ldots,l_{r}\right\}\subseteq V$ and the length of the path from $l_{i}$ to $v$ by $L_{i}$ for each $i\in\left\{1,\ldots,r\right\}$. We distinguish between even and odd lengths. Let $s\in \left\{0,\ldots,r\right\}$ be the number of even lengths and $t\in \left\{0,\ldots,r\right\}$ be the number of odd lengths such that $r=s+t$. Either one of the two sets of lengths is empty (i. e. for $s=0$ or $t=0$) or we have a set $\left\{L_{11},\ldots,L_{1s}\right\}$ of even lengths as well as a set $\left\{L_{21},\ldots,L_{2t}\right\}$ of odd lengths.
\label{Spid1}
\end{Sit}

Now, we can easily give an explicit formula for the cardinality of the biggest block of an independent partition of a spider.

\begin{Thm}
Given the conditions of Situation \ref{Spid1}.
\begin{enumerate}
\item[i)]
If all lengths are even, then 
\begin{eqnarray*}
M&=&\sum_{k=1}^{r}\frac{L_{k}}{2}
\end{eqnarray*}
holds.
\item[ii)]
If all lengths are odd, then we have
\begin{eqnarray*}
M&=&\sum_{k=1}^{r}\frac{L_{k}-1}{2}+1\;.
\end{eqnarray*}
\item[iii)]
If there are even and odd lengths, then 
\begin{eqnarray*}
M&=&\sum_{k=1}^{s}\frac{L_{1k}}{2}+\sum_{k=1}^{t}\frac{L_{2k}-1}{2}
\end{eqnarray*}
holds.
\end{enumerate}
\label{Spid2}
\end{Thm}

\textbf{Proof}:\\
\begin{itemize}
\item[i)]
For all $k\in\left\{1,\ldots,r\right\}$ let $P_{k}$ be the path from $l_{k}$ to $v$. By assumption, $2$ is a divisor of $L_{k}$. For $B\in\sigma$, $\sigma\in\Pi_{I}(S)$ with $|B|$ maximal we start to select vertices in $P_{k}$. Let $V(P_{k})=\left\{l_{k},v_{k1},v_{k2},\ldots,v_{kL_{k-2}},v\right\}$. According to Theorem \ref{MaxBlock1}, we start with $l_{k}$ because it is a leaf in $S$. Then $l_{k},v_{k2},v_{k4},\ldots,v_{kL_{k-2}}\in B$ and of course $|\left\{l_{k},v_{k2},v_{k4},\ldots,v_{kL_{k-2}}\right\}|=\frac{L_{k}}{2}$. Because we started with $l_{k}$ we have $v\notin B$ and hence we do not need to take account of $v$.
\item[ii)]
Here we apply the same method as abowe. Now, $L_{k}-1$ is even and we get the sequence $l_{k},v_{k2},v_{k4},\ldots,v_{kL_{k-3}}\in B$ with $|\left\{l_{k},v_{k2},v_{k4},\ldots,v_{kL_{k-3}}\right\}|=\frac{L_{k}-1}{2}$. But now $v\notin N(v_{kL_{k-3}})$ and thus we also choose $v\in B$. Because $v=\bigcap_{k=1}^{r}S_{k}$ we add exactly $1$.
\item[iii)]
We just combine the first two parts. Because there are even lengths we can not add $1$ for $v$.
\end{itemize}
q.e.d.\\

\end{document}